\documentclass[12pt,a4paper]{article}
\usepackage[latin2]{inputenc}
\usepackage{amsmath}
\usepackage{amsfonts}
\usepackage{amssymb}
\usepackage{graphicx}
\usepackage[left=2.00cm, right=2.00cm, top=2.50cm, bottom=2.50cm]{geometry}
\usepackage[T1]{fontenc}

\def\Q{{\mathbb Q}}
\def\Z{{\mathbb Z}}

\newtheorem{lemma}{Lemma}
\newtheorem{theorem}[lemma]{Theorem}

\title{
Calculating "small" solutions \\ of inhomogeneous relative Thue inequalities
}
\author{
Istv\'{a}n Ga\'{a}l\thanks{
        Research supported in part by the EFOP-3.6.1-16-2016-00022 project. The project is co-financed by the European Union and the European Social Fund.},\; \\
{\small University of Debrecen, Mathematical Institute} \\
{\small H--4002 Debrecen Pf.400., Hungary,} 
{\small e--mail: gaal.istvan@unideb.hu},
}

\begin{document}
\baselineskip=17pt

\maketitle
\thispagestyle{empty}

\renewcommand{\thefootnote}{\arabic{footnote}}
\setcounter{footnote}{0}


\begin{abstract}
Thue equations and their relative and inhomogeneous extensions are well known
in the literature. 
There exist methods, usually tedious methods, for the complete resolution of these equations. On the other hand our experiences show that such equations usually do 
not have extremely large solutions.
Therefore in several applications it is useful to have a fast algorithm
to calculate the "small" solutions of these equations. Under "small" solutions
we mean the solutions, say, with absolute values or sizes $\leq 10^{100}$.
Such algorithms were formerly constructed for Thue equations,
relative Thue equations. 
The relative and inhomogeneous 
Thue equations have applications in solving index form equations
and certain resultant form equations.
It is also known that certain "totally real" relative Thue equations can be reduced to
absolute Thue equations (equations over $\Z$).

As a common generalization of the above results,
in our paper we develop a fast algorithm for calculating "small" solutions 
(say with sizes $\leq 10^{100}$) of inhomogeneous relative Thue equations, more 
exactly of certain inequalities that generalize those equations. We shall show that 
in the "totally real" case these can similarly be reduced to absolute inhomogeneous
Thue inequalities. We also give an application to solving certain resultant 
equations in the relative case.
\end{abstract}

\newpage

\section{Introduction}

Throughout this paper we shall denote by $\Z_K$ the ring of integers of an
algebraic number field $K$. Further, $\overline{|\gamma|}$ will denote
the size of an algebraic number $\gamma$, that is the 
maximum of the absolute values of its conjugates.

\subsection{Thue equations, homogeneous, inhomogeneous, relative}
\label{ethue}

Let $F(x,y)\in\Z[x,y]$ be a homogenous form of degree $\ge 3$, irreducible over $\Q$,
and let $0\ne m\in\Z$. In 1909 A. Thue \cite{thue} proved that the equation
\begin{equation}
F(x,y)=m\;\; {\rm in}\;\; x,y\in\Z
\label{F}
\end{equation}
has only finitely many solutions. In 1968 A. Baker (see \cite{baker}) 
gave effective upper bounds for the solutions. These bounds have several improvements
but they all share the property that even in the simplest cases they are of magnitude
$\exp(10^{20})$, not allowing to calculate the solutions themselves. An algorithm to
reduce these bounds in order to explicitly determine the solutions of cubic Thue equations 
was first given by A. Baker and H. Davenport \cite{bada} which was extended to
arbitrary degrees by A. Peth\H o and R. Schulenberg \cite{psch}, and B. M. M. de Weger
\cite{deweger}. These algorithms give all solutions, however they require tedious computation.

It turned out that such equations usually only have a few small solutions 
(usually at most 3 digits) therefore a fast algorithm of A. Peth\H o \cite{pet} 
supplying "small" solutions turned to be very useful. 
Under "small" solutions we mean here solutions with
absolute values up to $10^{100}$--$10^{500}$, or even larger
(but definitely $\exp(10^{20})$ can not be reached). 
This algorithm produces within very short time 
all solutions that are applicable in practice but do not prove the non-existence 
of extremely large solutions (which, according to our experience do not exists
anyway).

If $\alpha$ is a root of $F(x,1)=0$ and $K=\Q(\alpha)$, then (\ref{F})
can be written as 
\begin{equation}
N_{K/Q}(x-\alpha y)=m\;\; {\rm in}\;\; x,y\in\Z.
\label{FN}
\end{equation}
As a generalization of these equations in 1974 V. G. Sprind\v zuk \cite{sprindzuk}
considered so called {\bf inhomogeneous Thue equations} of type
\begin{equation}
N_{K/Q}(x-\alpha y+\lambda)=m,
\label{FNL}
\end{equation}
in the variables $x,y\in\Z,\lambda\in\Z_K$, under the assumption
that $\lambda$ is a {\bf non-dominating variable}, that is
\[
\overline{|\lambda|}\leq (\max(|x|,|y|))^{1-\zeta},
\]
and $0<\zeta<1$ is a given.
V. G. Sprind\v zuk \cite{sprindzuk} gave upper bounds for the solutions, using Baker's 
method, which implies the finiteness of the number of solutions but do not enable one
to calculate the solutions.

Combining the method of \cite{sprindzuk} with \cite{psch}, \cite{deweger},
I. Ga\'al \cite{inhomThue} gave an algorithm to calculate all solutions of
inhomogeneous Thue equations of type (\ref{FNL}).

Note that inhomogeneous The equations gain important applications in solving
index form equations in sextic and octic fields with quadratic subfields, see
\cite{cycsex}, \cite{gpoctic}, \cite{book},
and also in solving certain resultant type equations, see Section \ref{rrx}.
In all these applications $\lambda$
is a fixed algebraic integer.

Thue equations were also extended to the so called {\bf relative case}, when
the form $F$ has coefficients in the ring of integers $\Z_M$ of an algebraic
number field $M$. These can be written as a relative norm form equation in two variables,
as
\begin{equation}
N_{K/M}(X-\alpha Y)=m\;\; {\rm in}\;\; X,Y\in\Z_M,
\label{FNM}
\end{equation}
where $\alpha$ is an algebraic number of degree $\geq 3$ over $M$, $K=M(\alpha)$
and $0\ne m\in\Z_M$.
Relative Thue equations were first considered by C. L. Siegel \cite{siegel}.
Upper bounds for the solutions were given by S. V. Kotov and V. G. Sprindzuk \cite{ks}.
An algorithm to calculate all solutions was constructed by
M. Pohst and I. Ga\'al \cite{gpthuerel}. A fast algorithm to compute "small" solutions of 
relative Thue equations was given by I. Ga\'al \cite{relthuesmall}. Further,
it turned out, that if $M$ is an imaginary quadratic field and all conjugates 
of $\alpha$ over $M$ are real, then the relative Thue equation can be reduced to 
absolute Thue equations (over $\Z$), cf. Section \ref{relabs}.

\subsection{Equations of type $F(x,y)=G(x,y)$}
\label{efg}

Thue equations (\ref{F}) were also extended to diophantine equations 
of type 
\begin{equation}
F(x,y)=G(x,y) \;\; {\rm in}\;\; x,y\in\Z,
\label{FG}
\end{equation}
where $F\in\Z[x,y]$ is as above, and $G\in\Z[x,y]$ 
with some restrictions on its degree, see the book of
T. N. Shorey and R. Tijdeman \cite{st}.
Fast algorithms for calculating "small" solutions of equations of type (\ref{FG})
were given by I. Ga\'al over $\Z$ , over imaginary quadratic fields, 
and also over arbitrary number fields (cf. 
\cite{FGZfelett}, \cite{fgim}, \cite{FGnumfield}, respectively).

\subsection{Application to resultant type equations}
\label{rrx}

An interesting application of inhomogeneous Thue equations is in solving 
{\bf resultant type equations}. 
Let $0\neq r\in\Z$. Determine polynomials $f,g\in\Z[x]$ with distinct roots
having 
\begin{equation}
{\rm Res}(f,g)=r.
\label{rs}
\end{equation}
Such diophantine problems were considered by W. M. Schmidt \cite{schres},
J. H. Evertse and K. Gy\H ory \cite{evgyres} and others.

It turned out, see I. Ga\'al \cite{resultantg=2}, that given $f$, the calculation of quadratic polynomials $g$ satisfying (\ref{rs}) leads to an inhomogeneous Thue 
equation, cf. Section \ref{rxrx}.

Given $f$, using the theory of unit equations 
I.Ga\'al and M.Pohst \cite{resultantnumfield}
gave an algorithm to determine all $g$ having roots in a fixed number field.

\section{Purpose}

In the present paper we are going to construct a fast algorithm to determine
"small" solutions of {\bf inhomogeneous relative Thue equations}. 
More exactly, as a common
generalization of the results in the previous Section we shall consider
inequalities.
These inequalities will have a polynomial on the left hand side 
having a structure similar to an inhomogeneous Thue form
\begin{equation}
N_{K/M}(X-\alpha Y+\lambda)
\label{xxyyll}
\end{equation}
($K=M(\alpha),[K:M]\geq 3, X,Y\in\Z_M, \lambda\in\Z_K$).
We consider the case when $\lambda$ is fixed, but this will be sufficient for
the applications in the resolution of the above mentioned index form equations
and resultant type equations. Further, since (\ref{xxyyll}) is the product of
inhomogeneous linear forms of type $X-\alpha_j Y +\lambda_j$, where $\alpha_j$, $\lambda_j$
are the relative conjugates of $\alpha$, $\lambda$, respectively, we shall
consider arbitrary products of inhomogeneous linear forms of this shape,
without requiring that $\alpha_j$, $\lambda_j$ be conjugates of certain 
algebraic numbers.

In order to cover equations of type (\ref{FG}),
on the right hand sides of these inequalities we shall have
some polynomial of $X,Y$. Since these polynomials
can be estimated from above by $\overline{|X|},\overline{|Y|}$, 
we use a power of
$Z=(\max(\overline{|X|},\overline{|Y|}))$ on the right hand side.

Our mail tool in the algorithm will be the application of the LLL
basis reduction algorithm, cf. \cite{lll}, \cite{dmv}.

\subsection{Reducing inhomogeneous relative Thue equations to absolute ones}
\label{relabs}

Using the ideas of I. Ga\'al, B. Jadrijevi\'c and L. Remete \cite{realthueimquadr}
in Section \ref{tot} we shall reduce totally real inhomogeneous relative Thue
equations to inhomogeneous Thue equations over $\Z$, that is, to the absolute case.
This can be done only in the special case, when all $\alpha_j$ are real, but in such
cases it is very useful for the resolution of the equations.

\subsection{Resultant type equations in the relative case}
\label{rxrx}

Let $M$ be an algebraic number field, let $f,g\in\Z_M[x]$ be polynomials,
irreducible over $M$, 
with roots $\alpha=\alpha_1,\alpha_2,\ldots,\alpha_n$, and
$\beta=\beta_1,\ldots,\beta_m$, respectively. 
The resultant of $f,g$ in this relative case
can be defined as
\[
{\rm Res}_M(f,g)=\prod_{i=1}^n\prod_{j=1}^m(\alpha_i-\beta_j).
\]
This can be written as
\[
\prod_{i=1}^n((\alpha_i-\beta_1)\cdots (\alpha_i-\beta_m))=
\prod_{i=1}^n g(\alpha_i)=N_{K/M}(g(\alpha)),
\]
where $K=M(\alpha)$.

Consider now the following problem. 
Let $f\in\Z_M[x]$ be a given polynomial,
having distinct roots $\alpha=\alpha_1,\alpha_2,\ldots,\alpha_n$.
Determine all monic quadratic polynomials $g(t)=t^2-Yt+X\in\Z_M[t]$,
with coefficients $X,Y\in\Z_M$, satisfying 
\begin{equation}
|{\rm Res}_M(f,g)|\leq c,
\label{fgx}
\end{equation}
with a given positive constant $c$.
By the above arguments (\ref{fgx}) can be written in the form
\[
|N_{K/M}(X-\alpha Y+\alpha^2)|\leq c \;\; {\rm in} \; X,Y\in\Z_M.
\]
Using our algorithm to determine "small" solutions of 
inhomogeneous relative Thue inequalities we shall be able to calculate the solutions
of (\ref{fgx}), with, say $\max(\overline{|X|},\overline{|Y|})\leq 10^{100}$.

\section{Inhomogeneous relative Thue inequalities}
\label{s1}

Let $M$ be a number field of degree $m$ with ring of integers $\Z_M$ 
and integral basis
$(\omega_1=1,\omega_2,\ldots,\omega_m)$. 
Assume that $X,Y$ are represented in the form 
\begin{equation}
X=x_1+x_2\omega_2+\ldots +x_m\omega_m, \; Y=y_1+y_2\omega_2+\ldots +y_m\omega_m,
\label{xy}
\end{equation}
with $x_j,y_j\in\Z\; (j=1,\ldots,m)$.
Set
\[
A=\max(\max_j |x_j|,\max_j |y_j|),
\]
and
\[
Z=(\max(\overline{|X|},\overline{|Y|})).
\]

Let $\alpha_1,\ldots,\alpha_n$ be given non-zero distinct complex numbers
and $\lambda_1,\ldots,\lambda_n$ any given complex numbers with
\begin{equation}
\max_j |\lambda_j|\leq c_{\lambda}.
\label{lambda}
\end{equation}
Let $c_0>0$ be a given positive constant.
As a generalization of the above equations consider the inequality
\begin{equation}
\left| \prod_{j=1}^n\left( X-\alpha_j Y +\lambda_j\right)\right|\leq c_0\cdot Z^k,
\;\;{\rm in}\;\; X,Y\in\Z_M,\;\; Z\leq Z_0,
\label{geneq}
\end{equation}
where $Z_0$ is huge given bound, say $10^{100}$.

Our arguments show, that the bound $Z\leq Z_0$ implies a bound $A_0$ for $A$ 
(of the same magnitude), which can be reduced in a numerical way to a much smaller value $A_R$ (say 100 or 1000),
such that all solutions of (\ref{geneq}) satisfy $A\leq A_R$
and these "tiny" possible solutions $X,Y$ with $A\leq A_R$ can be enumerated and tested.

In case $\alpha$ is an algebraic integer of degree $n$ over $M$,
with relative conjugates $\alpha_1,\ldots,\alpha_n$ over $M$
and $\lambda$ is an algebraic integer in
$K=M(\alpha)$ with relative conjugates of $\lambda_1,\ldots,\lambda_n$ over $M$, then
inequality (\ref{geneq}) can be written as
\begin{equation}
|N_{K/M}(X-\alpha Y+\lambda)|\leq c_0\cdot Z^k,
\label{norm}
\end{equation}
which can be considered as an inhomogeneous relative Thue inequality,
analogues to the one considered by V. G. \v Sprindzuk.

\subsection{Elementary estimates}
\label{s2}

The basic tool in solving various types of Thue equations is that if the
solutions are not very small, then a linear factor of the left hand side
must be very small. This is what we shall show in this section using estimates
that are more or less standard in the theory of Thue equations.

Let $0<\varepsilon_i<1\; (1\leq i\leq n)$. These we can choose appropriately, see our remarks in the proof of 
Lemma \ref{lem1}.
For $1\leq i,j\leq n,i\ne j$ set
\[
c_{1ij}=|\alpha_j-\alpha_i|\min\left(1,\frac{1}{|\alpha_i|}\right),\;\;
c_{2ij}=c_0^{1/n}\max\left(1,\frac{|\alpha_j|}{|\alpha_i|}\right),\;\;
c_{3ij}=\max\left(|\lambda_j-\lambda_i|,
\frac{|\alpha_i\lambda_j-\alpha_j\lambda_i|}{|\alpha_i|}\right),
\]
\begin{equation}
c_{4i}=\max_{j\neq i}\left(\max\left(
\left(\frac{2c_{2ij}}{\varepsilon_i c_{1ij}}\right)^{\frac{n}{n-k}},
\frac{2c_{3ij}}{(1-\varepsilon_i) c_{1ij}}
\right)\right),
\label{c4}
\end{equation}
\[
c_{5i}=\frac{2^{n-1}c_0}{\prod_{j\ne i}c_{1ij}}.
\]

Let $X,Y\in \Z_M$ be a solution of (\ref{geneq}). 
Let $i\; (1\leq i\leq n)$ be the index with
\begin{equation}
|X-\alpha_i Y +\lambda_i|=\min_j |X-\alpha_j Y +\lambda_j|.
\label{ii}
\end{equation}
This $i$ is different for different solutions. Therefore we have to consider
all possible values of $i$ to find all solutions of (\ref{geneq}), 
that is, the following procedure
must be performed for each $i$.

\begin{lemma}
Assume that $X,Y\in \Z_M$ is a solution of (\ref{geneq}) with (\ref{ii}).
If 
\begin{equation}
Z\geq c_{4i},
\label{ZZ}
\end{equation}
then 
\begin{equation}
|\beta_i|\leq c_{5i}\cdot Z^{k-n-1}.
\label{betaismall}
\end{equation}
\label{lem1}
\end{lemma}

\noindent
{\bf Proof of Lemma \ref{lem1}}.
Set $\beta_j=X-\alpha_j Y +\lambda_j\; (1\leq j\leq n)$. Then (\ref{ii}) yields
\begin{equation}
|\beta_i|=\min_j |\beta_j|.
\label{betai}
\end{equation}
Inequality (\ref{geneq}) implies
\begin{equation}
|\beta_i|\leq c_0^{1/n}\cdot Z^{1/n}.
\label{betai2}
\end{equation}
On the other hand, for $j\ne i\; (1\leq j\leq n)$ we have
\[
|\beta_j|\geq |\beta_j-\beta_i|-|\beta_i|\geq |\alpha_j-\alpha_i|Y
-|\lambda_j-\lambda_i|- c_0^{1/n} Z^{k/n}.
\]
Further,
\[
\alpha_i\beta_j=(\alpha_i-\alpha_j)X+\alpha_i\lambda_j-\alpha_j\lambda_i+\alpha_j\beta_i,
\]
whence
\[
|\beta_j|\geq\frac{|\alpha_i-\alpha_j|}{|\alpha_i|}X
-\frac{|\alpha_i\lambda_j-\alpha_j\lambda_i|}{|\alpha_i|}
-\frac{|\alpha_j|}{|\alpha_i|}c_0^{1/n} Z^{k/n}.
\]
The above inequalities imply
\begin{equation}
|\beta_j|\geq c_{1ij}Z-c_{2ij}Z^{k/n}-c_{3ij}.
\label{betaj}
\end{equation}
Our estimate (\ref{betaj}) implies, that if $Z$ is large enough, then 
\begin{equation}
|\beta_j|\geq \frac{c_{1ij}}{2}Z.
\label{betaj2}
\end{equation}
This is satisfied if
\[
c_{2ij}Z^{k/n}+c_{3ij}\leq \frac{c_{1ij}}{2}Z.
\]
We have to cover both summands by a part of $c_{1ij}Z/2$. Here 
we make use of $0<\varepsilon_i <1$ and assume that
\[
c_{2ij}Z^{k/n}\leq  \varepsilon_i \frac{c_{1ij}}{2}Z,\;\;\; 
c_{3ij}\leq (1-\varepsilon_i)\frac{c_{1ij}}{2}Z,
\]
which are satisfied in view of (\ref{ZZ}).
An optimal choice of $\varepsilon_i$ is useful to keep the constant $c_{4i}$ 
in (\ref{c4}), (\ref{ZZ})
as small as possible. 
In view of (\ref{betaj2}) our inequality (\ref{geneq}) implies (\ref{betaismall}).
\hfill $\Box$

\subsection{Comparing $X,Y$ with $x_j,y_j$}
\label{s3}

Denote by $\omega_k^{(j)}$ the conjugates of $\omega_k$,
($1\leq j,k \leq m$).
Set
\[
S=\left(
\begin{array}{cccc}
1&\omega_2^{(1)}&\hdots&\omega_m^{(1)}\\
\vdots&\vdots&&\vdots\\
1&\omega_2^{(m)}&\hdots&\omega_m^{(m)}
\end{array}
\right).
\]
Denote by $c_6$ the row-norm of the matrix $S$, that is the maximum sum
of the absolute values of the entries in its rows.
Let $c_7$ be the row norm of the matrix $S^{-1}$.
For $1\leq i\leq n$ let
\[
c_{8i}=\frac{c_{4i}}{c_7},\;\; c_{9i}=c_{5i}\cdot c_7^{n+1-k}.
\]

\begin{lemma}
Assume that $X,Y\in \Z_M$ in the representation(\ref{xy}) 
is a solution of (\ref{geneq}) with (\ref{ii}). 
We have 
\begin{equation}
Z\leq c_6\cdot A.
\label{ZA}
\end{equation}
If 
\begin{equation}
A\geq c_{8i}.
\label{A>}
\end{equation}
then
\begin{equation}
|\beta_i|\leq c_{9i} A^{k-n-1}.
\label{bA}
\end{equation}
Moreover, if $Z< Z_0$ then
\begin{equation}
A\leq A_0= c_7 Z_0.
\label{A<}
\end{equation}
\label{lem2}
\end{lemma}

\noindent
{\bf Proof of Lemma \ref{lem2}}.
We have
\begin{equation}
\left(
\begin{array}{c}
X^{(1)}\\ \vdots \\ X^{(m)}
\end{array}
\right)
=
S\cdot
\left(
\begin{array}{c}
x_1\\ \vdots \\ x_m
\end{array}
\right).
\label{S}
\end{equation}
Obviously
\[
\overline{|X|}\leq c_6\cdot \max_j |x_j|,
\]
and similarly
\[
\overline{|Y|}\leq c_6\cdot \max_j |y_j|,
\]
whence (\ref{ZA}) follows.
Further, 
\[
\left(
\begin{array}{c}
x_1\\ \vdots \\ x_m
\end{array}
\right)
=
S^{-1}
\cdot
\left(
\begin{array}{c}
X^{(1)}\\ \vdots \\ X^{(m)}
\end{array}
\right),
\]
which implies 
\[
\max_j(|x_j|)\leq c_7 \overline{|X|},
\]
and similary 
\[
\max_j(|y_j|)\leq c_7 \overline{|Y|}.
\]
Therefore we obtain
\begin{equation}
A\leq c_7 Z.
\label{XYZ}
\end{equation}
Therefore $Z< Z_0$ implies (\ref{A<}).

(\ref{XYZ}) and (\ref{A>}) imply (\ref{ZZ}), whence,
in view of Lemma \ref{lem1} we obtain (\ref{betaismall}).
Finally, by (\ref{XYZ}) inequality (\ref{betaismall}) implies (\ref{bA}).
$\Box$

\subsection{Reducing the bounds}
\label{s4}

In view of our Lemma 1 and Lemma 2 we are going to calculate 
the solutions $x_1,\ldots,x_m,y_1,\ldots,y_m\in \Z$ of (\ref{bA}) satisfying 
(\ref{A<}).

Let $H$ be a large constant to be specified later, let $i$ be the index with
(\ref{ii}).
In view of Lemma \ref{lem2} we consider the inequality 
\begin{equation}
|x_1+\omega_2 x_2+\ldots + \omega_m x_m 
-\alpha_i y_1-\alpha_i \omega_2 y_2-\ldots - \alpha_i\omega_m y_m+\lambda_i | \leq 
c_{9i}\cdot A^{k-n-1}.
\label{redineq}
\end{equation}

Consider now the lattice $\cal{L}$ generated by the columns of the matrix
{\small
\[
\cal{L}=
\left(
\begin{array}{ccccccccccc}
1&0&\ldots&0&\vline&0&0&\ldots &0&\vline &0\\
0&1&\ldots&0&\vline&0&0&\ldots &0&\vline &0\\
\vdots&\vdots&\ddots&\vdots&\vline&\vdots&\vdots&\ddots&\vdots&\vline&\vdots\\
0&0&\ldots&1&\vline&0&0&\ldots &0&\vline &0\\
\hline
0&0&\ldots&0&\vline&1&0&\ldots &0&\vline &0\\
0&0&\ldots&0&\vline&0&1&\ldots &0&\vline &0\\
\vdots&\vdots&\ddots&\vdots&\vline&\vdots&\vdots&\ddots&\vdots&\vline&\vdots\\
0&0&\ldots&0&\vline&0&0&\ldots &1&\vline &0\\
\hline
0&0&\ldots&0&\vline&0&0&\ldots &0&\vline &1\\
\hline
H&H{\rm Re}(\omega_2)&\ldots & H{\rm Re}(\omega_m)&\vline&
H{\rm Re}(\alpha_i)&H{\rm Re}(\alpha_i\omega_2)&\ldots & H{\rm Re}(\alpha_i\omega_m)&\vline&H{\rm Re}(\lambda_i) \\
H&H{\rm Im}(\omega_2)&\ldots & H{\rm Im}(\omega_m)&\vline&
H{\rm Im}(\alpha_i)&H{\rm Im}(\alpha_i\omega_2)&\ldots & H{\rm Im}(\alpha_i\omega_m)&\vline&H{\rm Im}(\lambda_i)\\
\end{array}
\right).
\]
}
In the totally real case we may omit the last row.
Denote by $b_1$ the first vector of an LLL-reduced basis of the lattice ${\cal L}$
(cf. \cite{lll},\cite{dmv}).

\begin{theorem}
\label{redlemma}
Assume that $X,Y\in \Z_M$ in the representation(\ref{xy}) 
is a solution of (\ref{geneq}) with (\ref{ii}). 
Assume that 
\begin{equation}
A=\max(\max_j|x_j|,\max_j|y_.|)\geq \max(c_{8,i}, 2c_{\lambda}).
\label{AAA}
\end{equation}
If $A\leq A_0$ for some positive constant $A_0$ and $H$ is choosen so large that
for the first vector $b_1$ of the LLL reduced basis of $\cal L$ we have
\begin{equation}
|b_1|\geq \sqrt{(2m+1)2^{2m}} \cdot A_0,
\label{b1}
\end{equation}
then
\begin{equation}
A\leq \left( \frac{2c_{8i} H}{\sqrt{3} A_0}\right)^{\frac{1}{n-k-1}}.
\label{rrr}
\end{equation}
\label{th3}
\end{theorem}

\noindent
{\bf Proof of Theorem \ref{th3}}.\\ 
We follow the arguments of \cite{relthuesmall} (see also \cite{book}).
Denote by $l_0$ the shortest non-zero vector in the lattice $\cal L$.
Denote $l_1$ the shortest non-zero vector in the lattice $\cal L$ having last
coordinate 1.
Assume that the vector $l$ is a linear combination of the
lattice vectors with coefficients $x_1,\ldots,x_m,y_1,\ldots,y_m,1$,
respectively. 

Using the inequalities of \cite{dmv} we have $|b_1|^2\leq 2^{2m} |l_0|^2$. 
Obviously $|l_0|\leq |l_1|\leq |l|$.
The first $2m$ components of $l$ are in absolute value $\leq A_0$,
for the last 2 components (\ref{redineq}) is satisfied. 
(Observe that the last two components of $l$ are 
the real and imaginary parts of the left hand side of (\ref{redineq})).

Hence we obtain
\begin{eqnarray*}
&&(2m+1)A_0^2 = 2^{-2m}\left( (2m+1)\cdot 2^{2m} A_0^2 \right)  \\
&&\leq 2^{-2m}|b_1|^2\leq |l_0|^2\leq |l_1|^2\leq |l|^2
\leq 2m \cdot A_0^2+c_{\lambda}^2+H^2 c_{8i}^2 A^{2(k+1-n)},
\end{eqnarray*}
(cf. \ref{lambda}) whence 
\[
A_0^2\leq c_{\lambda}^2+c_{8i}^2 H^2 A^{2(k+1-n)}.
\]
In view of (\ref{AAA}) we obtain 
\[
\frac{3}{4}A_0^2\leq c_{8i}^2 H^2 A^{2(k+1-n)},
\]
which implies the assertion (\ref{rrr}).
\hfill $\Box$\\

\subsection{How to apply Lemma \ref{redlemma} ?}
\label{s5}

If $Z\leq Z_0$, then by (\ref{A<}) we have $A\leq A_0$ with $A_0=c_7Z_0$.
We make $H$ so large, that (\ref{b1}) is satisfied. Then Lemma \ref{redlemma} implies
a new bound (\ref{rrr}) for $A$. 

In the next step this new bound for $A$ will play the role of $A_0$.
If $H$ is large enough and (\ref{b1}) is satisfied, then 
by (\ref{rrr}) we get a new improved bound for $A$.
We continue this process, until the new bound is less than the previous one.

Usually we make several tests to choose $H$ appropriately. 
According to our experiences, the order of magnitude of $H$
must be about $A_0^{2m+1}$. In the first reduction steps $A_0$ and $H$
are much larger then the other constants in the formula (\ref{rrr}).
Therefore the constant on the right side of (\ref{rrr}) is of magnitude
$A_0^{2m/(n-k-1)}$. Hence the reduction is only efficient, if
\[
2m<n-k-1,
\]
that is
\[
n>2m+k+1.
\]
Our examples give a good insight how the reduction procedure works.

\subsection{The "tiny" possible values of the variables}

It is important to remark, that we 
have to deal separately with those $x_1,\ldots,x_m,
y_1,\ldots,y_m$ not satisfying (\ref{AAA}).
As we shall see later in our examples, a sensitive part of the algorithm is to 
enumerate these "tiny" possible values of the variables.

\section{Totally real inhomogeneous relative Thue inequalities over 
imaginary quadratic fields}
\label{tot}

We have seen in \cite{realthueimquadr} that if the 
ground field $M$ is a complex quadratic field and the coefficients of the 
relative Thue equation are real, then the relative Thue inequality can be
reduced to absolute Thue inequalities (over $\Q$). This is what we extend to the 
inhomogeneous case in this section. Our arguments applied to 
\cite{realthueimquadr} may improve some of those estimates.

Let $d>1$ be a square-free positive integer, and set $M=\Q(i\sqrt{d})$. 

\noindent
If $d\equiv 3\; (\bmod \; 4)$, then $X,Y\in\Z_M$ can be written as
\[
X=x_1+x_2\frac{1+i\sqrt{d}}{2}=\frac{(2x_1+x_2)+x_2i\sqrt{d}}{2}, 
\]
\[
Y=y_1+y_2\frac{1+i\sqrt{d}}{2}=\frac{(2y_1+y_2)+y_2i\sqrt{d}}{2},
\]
with $x_1,x_2,y_1,y_2\in\Z$.

\noindent
If $d\equiv 1,2\; (\bmod \; 4)$, then $X,Y\in\Z_M$ can be written as
\[
X=x_1+x_2i\sqrt{d},\;\; Y=y_1+y_2i\sqrt{d},
\]
with $x_1,x_2,y_1,y_2\in\Z$.

Let $\alpha_1,\ldots,\alpha_n$ be distinct non-zero real numbers,
and let $\lambda_j=\lambda_{j1}+i\lambda_{j2}$  where
$\lambda_{j1},\lambda_{j2}$ are real numbers ($1\leq j\leq n$)
and $i=\sqrt{-1}$.
Let $c_0$ be a positive constant.

\begin{theorem}
Assume
\begin{equation}
\left|\prod_{j=1}^n(X-\alpha_j Y+\lambda_j)\right|\leq c_0.
\label{re}
\end{equation}
In case $d\equiv 3\; (\bmod \; 4)$, we have
\begin{equation}
\left|\prod_{j=1}^n((2x_1+x_2)-\alpha_j (2y_1+y_2)+2\lambda_{j1})\right|\leq 2^n c_0,
\label{re31}
\end{equation}
and 
\begin{equation}
\left|\prod_{j=1}^n(x_2-\alpha_j y_2+\frac{2}{\sqrt{d}}\lambda_{j2})\right|
\leq \frac{2^nc_0}{\sqrt{d}^n}.
\label{re32}
\end{equation}
In case $d\equiv 1,2\; (\bmod \; 4)$, we have
\begin{equation}
\left|\prod_{j=1}^n(x_1-\alpha_j y_1+\lambda_{j1})\right|\leq c_0,
\label{re11}
\end{equation}
and
\begin{equation}
\left|\prod_{j=1}^n(x_2-\alpha_j y_2+\frac{2}{\sqrt{d}}\lambda_{j1})\right|
\leq \frac{c_0}{\sqrt{d}^n}.
\label{re12}
\end{equation}
\label{lemim}
\end{theorem}

\vspace{0.5cm}

Observe, that the inequalities 
(\ref{re31}), (\ref{re32}), (\ref{re11}), (\ref{re12}) are all over $\Z$
(the absolute case).

\vspace{0.5cm}

\noindent
{\bf Proof of Theorem \ref{lemim}}. 
Let
\[
\beta_j=X-\alpha_j Y+\lambda_j,
\]
for $1\leq j\leq n$. The proof is based on the simple observation that
\[
|{\rm Re}(\beta_j)|\leq |\beta_j|,\;\; |{\rm Im}(\beta_j)|\leq |\beta_j|,
\]
therefore
\[
\prod_{j=1}^n|{\rm Re}(\beta_j)|\leq \prod_{j=1}^n|\beta_j|\leq c_0,
\]
and
\[
\prod_{j=1}^n|{\rm Im}(\beta_j)|\leq \prod_{j=1}^n|\beta_j|\leq c_0.
\]
\hfill$\Box$

The "small" solutions of inequalities (\ref{re31}), (\ref{re32}),  (\ref{re11}), (\ref{re12})
can be calculated by applying the procedure of the previous sections with $m=1$.
A typical application will be given in Section \ref{rreess}.

\section{Examples}

We illustrate our procedures by two detailed examples.
In our examples we have $k=0$ (cf. (\ref{geneq})) in order to make the
reduction procedures more efficient.

\subsection{An inhomogeneous relative Thue equation of degree 9\\ over a real 
quadratic field}

Let $\alpha=\alpha_1,\alpha_2,\ldots,\alpha_9$ be the roots of 
\[
f(x)=x^9-9x^7+24x^5-2x^4-20x^3+3x^2+5x-1.
\]
Note that all roots are real. This polynomial was taken from \cite{voight}.
Let $\lambda_i=\alpha_i^2+2\alpha_i\; (1\leq i\leq 9)$.
Let $M=\Q(\sqrt{2})$ and consider the inhomogeneous
relative Thue inequality
\[
\left|\prod_{j=1}^9(X-\alpha_i Y+\lambda_i)\right|\leq 10 \;\; {\rm in}\;\; X,Y\in\Z_M.
\]
Setting $K=M(\alpha)$ and $\lambda=\alpha^2+2\alpha$ this can be written as
\[
\left|N_{K/M}(X-\alpha Y +\lambda)\right|\leq 10 \;\; {\rm in}\;\; X,Y\in\Z_M.
\]
We are going to determine all solutions $X,Y\in\Z_M$ with 
\[
\max(\overline{|X|},\overline{|Y|})\leq 10^{100}.
\]
This bound gives 
\[
A=\max(|x_1|,|x_2|,|y_1|,|y_2|)\leq 10^{100}.
\]
For all $i$ we took $\varepsilon_i=0.5$. In view of (\ref{AAA}) our estimates 
are valid for $A\geq 56$.

The following table details the reduction procedure. 
In the columns of the table we display 
the number of step, 
the initial bound $A_0$ for $A$,
the lower bound for $||b_1||$, 
the value of $H$ used in that step,
the precision (number of digits) used in that step
and the reduced bound.
\[
\begin{array}{|c|c|c|c|c|c|c|}
\hline
{\rm step}& A_0  & ||b_1||\geq    &    H   &  {\rm Digits}   &  {\rm new}\; A_0  \\ \hline
1.&10^{100} & 4.4721\cdot 10^{100} & 10^{505}  &  600 & 9.0912\cdot 10^{50}   \\ \hline 
2.&9.0912\cdot 10^{50} & 4.0657\cdot 10^{51} & 10^{260}  & 350 &2.9093\cdot 10^{26}\\ \hline 
3.&2.9093\cdot 10^{26} & 1.3010\cdot 10^{27} & 10^{137}  & 250 &1.4146\cdot 10^{14}\\ \hline 
4.&1.4146\cdot 10^{14} & 6.3262\cdot 10^{14} & 10^{77}  & 150 &1.5480\cdot 10^{8}\\ \hline 
5.&1.5480\cdot 10^{8} & 6.9228\cdot 10^{8} & 10^{46}  & 100 &1.1478\cdot 10^{5}\\ \hline 
6.&1.1478\cdot 10^{5} & 5.1331\cdot 10^{5} & 10^{30}  & 100 &2825 \\ \hline 
7.& 2825  & 12633.7840 & 10^{22}  & 50 & 449 \\ \hline 
8.& 449  & 2007.9890 & 10^{18}  & 50 & 178 \\ \hline 
9.& 178  & 796.0401 & 10^{15}  & 50 & 84 \\ \hline 
10.& 84  & 375.6594 & 10^{14}  & 50 & 69 \\ \hline 
11.& 69  & 308.5773 & 4\cdot 10^{13}  & 50 & 63 \\ \hline 
12.& 63  & 281.7445 & 3\cdot 10^{13}  & 50 & 62 \\ \hline 
13.& 62  & 277.2724 & 2.9\cdot 10^{13}  & 50 & 61 \\ \hline 
\end{array}
\]
The possible $x_1,x_2,y_1,y_2$ with absolute values $\leq 61$ were enumerated.
We found 138 solutions, for all of them we have 
$|x_1|,|x_2|,|y_1|,|y_2|\leq 4$.

\subsection{Application to relative resultant inequalities}
\label{rreess}

Here we give an application of the results in Section \ref{tot}
to solving a relative resultant inequality (cf. Section \ref{rxrx}).

Set $M=\Q(i\sqrt{2})$. Let 
\[
f(t)=t^5-t^4-4t^3+3t^2+3t-1.
\]
This totally real quintic polynomial was taken from \cite{diaz}.
We shall consider this polynomial as $f(t)\in\Z_M[t]$.
Our task is now to determine all monic quadratic polynomials 
$g(t)=t^2-Yt+X\in\Z_M[t]$ with 
\begin{equation}
|{\rm Res}_M(f,g)|\leq 25,
\label{fg}
\end{equation}
and
\begin{equation}
\max(\overline{|X|},\overline{|Y|})<10^{100}.
\label{XXYY}
\end{equation}
Denote by $\alpha=\alpha_1,\alpha_2,\ldots,\alpha_5$ the roots of $f(t)\in\Z_M[t]$,
and by $\beta_1,\beta_2$ the roots of $g(t)\in \Z_M[t]$.

As we have seen in Section \ref{rxrx}, inequality (\ref{fg}) can be written as
\begin{equation}
|N_{K/M}(X-\alpha Y+\alpha^2)|\leq 25 \;\; {\rm in} \; X,Y\in\Z_M,
\label{rree}
\end{equation}
where $K=\Q(\alpha,i\sqrt{2})$. 
Set $X=x_1+i\sqrt{2}x_2,Y=y_1+i\sqrt{2}y_2$. 
Let $L=\Q(\alpha)$. In view of Theorem \ref{lemim} inequality (\ref{rree}) implies
\begin{equation}
|N_{L/Q}(x_1-\alpha y_1+\alpha^2)|\leq 25,\;\; {\rm in} \;x_1,y_1\in\Z,
\label{inhom}
\end{equation}
and
\begin{equation}
|N_{L/Q}(x_2-\alpha y_2)|\leq \frac{25}{(\sqrt{2})^5},\;\; {\rm in} \;x_2,y_2\in\Z.
\label{hom}
\end{equation}
In view of the bound given for $X,Y$,
we have to find the solutions of these inequalities with
\begin{equation}
\max(|x_1|,|x_2|,|y_1|,|y_2|)\leq 10^{100}.
\label{korl}
\end{equation}
The second inequality is a simple Thue inequality, the "small" solutions of which can be
easily calculated by using \cite{pet}.

The first inequality can be solved by using the algorithm of Section
\ref{s1} with $m=1$.
We detail here the reduction procedure, similarly as in the previous example.

\[
\begin{array}{|c|c|c|c|c|c|c|}
\hline
{\rm step}& A_0  & ||b_1||\geq    &    H   &  {\rm Digits}   &  {\rm new}\; A_0  \\ \hline
1.&10^{100} & 4.8989\cdot 10^{100} & 10^{303}  &  400 & 2.8230\cdot 10^{51}   \\ \hline 
2.&2.8230\cdot 10^{51} & 1.3829\cdot 10^{52} & 10^{160}  & 250 &6.8871\cdot 10^{27}\\ \hline 
3.&6.8871\cdot 10^{27} & 3.3739\cdot 10^{28} & 10^{87}  & 150 &3.0989\cdot 10^{15}\\ \hline 
4.&3.0989\cdot 10^{15} & 1.5181\cdot 10^{16} & 10^{50}  & 100 &2.1277\cdot 10^{9}\\ \hline 
5.&2.1277\cdot 10^{9} & 1.0423\cdot 10^{10} & 10^{31}  & 100 &1.3144\cdot 10^{6}\\ \hline 
6.&1.3144\cdot 10^{6} & 6.4392\cdot 10^{6} & 2\cdot 10^{21}  & 100 &31354 \\ \hline 
7.& 31354  & 1.5360\cdot 10^5 & 20\cdot 10^{15}  & 50 & 4486 \\ \hline 
8.& 4486  & 21976.8219 & 4\cdot 10^{13}  & 50 & 1542 \\ \hline 
9.& 1542  & 7554.2263 & 17\cdot 10^{11}  & 50 & 914 \\ \hline 
10.& 914  & 4477.6672 & 4\cdot 10^{11}  & 50 & 726 \\ \hline 
11.& 726  & 3556.6591 & 10^{11}  & 50 & 543 \\ \hline 
12.& 543  & 2660.1458 &  6\cdot 10^{10}  & 50 & 514 \\ \hline 
13.&  514 & 2518.0754 & 5\cdot 10^{10}  & 50 & 498 \\ \hline 
14.&  498 & 2439.6917 & 4.5\cdot 10^{10}  & 50 & 489 \\ \hline 
15.&  489 & 2395.6009& 4.2\cdot 10^{10}  & 50 & 483 \\ \hline 
16.&  483 & 2366.2070 & 4.1\cdot 10^{10}  & 50 & 481 \\ \hline 
\end{array}
\]
The reduced bound 481 is larger then (\ref{AAA}), which gives 45.
We have enumerated all $x_1,y_1$ with absolute values $\leq 481$. For all solutions 
of (\ref{inhom}) we obtained $|x_1|,|y_1|\leq 4$. Using the method of \cite{pet}
we calculated the solutions of (\ref{hom}) with (\ref{korl}) and we obtained
a few solutions with $|x_2|,|y_2|\leq 2$. From the above values of 
$x_1,x_2,y_1,y_2$ we obtained 39 solutions of (\ref{rree}), which allowed
to create all polynomials $g(t)=t^2-(y_1+i\sqrt{2}y_2)t+(x_1+i\sqrt{2}x_2)\in\Z_M[t]$
satisfying (\ref{fg}) with (\ref{XXYY}).

\subsection{Computational aspects}

All calculations were performed in Maple \cite{maple} and were executed 
on an average laptop running under Windows. The CPU time took some seconds,
including also the reduction steps, except for the enumeration of "tiny" 
values after the reduction procedures which took 1-2 minutes.

\end{document}